\documentclass[12pt,twoside]{amsart}
\usepackage{amsmath}
\usepackage{amsfonts}
\usepackage{amssymb}
\usepackage{amsthm}
\usepackage[english]{babel}

\newtheorem{theo}{Theorem}
\newtheorem{lem}[theo]{Lemma}
\newtheorem{prop}[theo]{Proposition}
\newtheorem{cor}[theo]{Corollary}

\newcommand{\E}{\ensuremath{\mathbb {E}}}

\newcommand{\A}{{\mathcal A}}
\newcommand{\F}{{\mathcal F}}

\def\E{\mathbb{E}}

\def\F{{\mathcal F}}

\begin{document}

\title[Quenched invariance principle]{ A quenched invariance principle for stationary processes }

\medskip

 \author{ Christophe Cuny}
\curraddr{Laboratoire MAS, Ecole Centrale de Paris, Grande Voie des Vignes, 92295 Chatenay-Malabry cedex
FRANCE}
\email{christophe.cuny@ecp.fr}

\author{Dalibor Voln\'y}
\curraddr{D\'epartement de Math\'ematiques, Universit\'e de Rouen,
76801 Saint-Etienne du Rouvray, FRANCE}
\email{dalibor.volny@univ-rouen.fr}
\thanks{}



\medskip

\begin{abstract}
In this note, we prove a conditionally centered version of the quenched weak invariance principle
under the Hannan condition, for stationary processes. In the course, we obtain a (new) construction of the fact that any
stationary process may be seen as a functional of a Markov chain.
\end{abstract}

\maketitle

\section{Introduction and results}

Let $(X,{\mathcal A},\mu)$ be a probability space and $\theta$ be an
invertible bimeasurable
transformation of $X$, preserving $\mu$, and assume that $\theta$ is ergodic. Let ${\mathcal F}_0$ be a
sub-$\sigma$-algebra
of ${\mathcal A}$ such that ${\mathcal F}_0\subset \theta^{-1}(
{\mathcal F}_0)$. Define a filtration  $({\mathcal F}_n)_{n\in {\bf Z}}$,
by ${\mathcal F}_n=\theta^{-n}{\mathcal F}_0$ and denote 
$\F_{-\infty}=\cap_{n\in \mathbb{Z}}\F_n$. For every $n\in\bar {\mathbb Z}$,
we denote by $\E_n$ the conditional expectation with respect to
${\mathcal F}_n$ and we define the projection ${P}_n:=
\E_n-\E_{n-1}$.

\medskip

Let $f$ be ${\mathcal F}_0$-measurable. We want to study the stationary process
$(f\circ \theta^n)_{n\in {\bf Z}}$.

\medskip

Let $\mu(\cdot,\cdot)$ denote a regular conditional probability on $\A$ given $\F_0$, see \cite{B} p. 358-364, and
for every $x\in X$, write
$\mu_x:= \mu(x,\cdot)$. Thus, for every $x\in X $, $\mu_x$ is a probability
measure on $\A$, and for every $A\in \A$, $\mu(\cdot,A)$ is a version of $\mu(A|\F_0)$.

\medskip

\medskip

We say that the process $(f\circ \theta^i)$ satisfies the 
\emph{Hannan Condition} if
\begin{equation}\label{Hannan}
  \sum_{i=0}^\infty \|P_if\|_2 = \sum_{i=0}^\infty \|P_0(f\circ \theta^i)\|_2 < \infty.
\end{equation}
\medskip

If $\E_{-\infty}(f)=0$, the Hannan Condition, introduced by E.J. Hannan in \cite{H}, guarantees the CLT and the weak invariance 
principle (WIP). The condition
has been shown to be very useful in applications (cf. \cite{DMV}, also for the WIP). 
In general, as shown in \cite{D}, the Hannan Condition is independent  of the so-called Dedecker-Rio and Maxwell-Woodroofe 
conditions, that are also sufficient for the WIP.

\medskip

Let us denote $S_n=S_n(f) = \sum_{i=1}^n f\circ \theta^i$.
It was shown in \cite{VW1} that  the
CLT is not quenched (i.e. $S_n$ does not satisfy the CLT under $\mu_x$ for $\mu$-almost every $x\in X$) 
but it follows from \cite{VW2} or \cite{CP} that
$\bar S_n=S_n - \E_0(S_n)$
satisfies the quenched CLT. Here we prove that $ \bar S_n$ satisfies  the quenched WIP  as well.
\medskip

For every $t\in [0,1]$, write $S_n(t)=S_{[nt]}+(nt-[nt])f\circ 
\theta^{[nt]+1}$ and $\bar S_n(t)=S_n(t)-\E_0(S_n(t))$. Our main result is the following.

\begin{theo}\label{theo}
Let $f\in L^2(X,\A,\mu)$ satisfy the Hannan condition. Then
there exists a martingale $(M_n)_n$ with stationary ergodic increments such that
\begin{equation}\label{strest}
\E_0(\max_{1\le n\le N}(\bar S_n-M_n)^2)= o(N) \qquad \mbox{$\mu$-a.s.}
\end{equation}
In particular, $\sigma^2:=\lim_n \E(S_n^2)/n$ exists and  for $\mu$-almost every $x\in X$, for every bounded continuous function $\varphi$ on 
$C([0,1], \|\cdot \|_\infty)$, we have
\begin{equation*}
\int_X \varphi(S_n(t )/\sqrt n)d\mu_x \underset{n\rightarrow +\infty}
\longrightarrow \E(\varphi(\sigma W_t)),
\end{equation*}
where $(W_t)_{0\le t\le 1}$ stands for a standard brownian motion.
\end{theo}

\begin{cor}\label{cor}
Let $f\in L^2(X,\A,\mu)$ be such that
\begin{equation}\label{SMW}
\sum_{n\ge 1} \frac{\|\E_0(X_n)\|_2}{\sqrt n}<\infty
\end{equation}
Then, \eqref{Hannan} holds and the conclusion of Theorem \ref{theo} 
is true with $S_n(t)$ in place of $\bar S_n(t)$.
\end{cor}



\section{Proof of the results}

To avoid technical difficulties (and since it is also convenient for the next 
section) we assume that $X$ is a Polish space and that $\A$ is the 
$\sigma$-algebra of its Borel sets. It is known (see for instance 
Neveu \cite[Proposition V.4.3]{Neveu}) 
that in this case there exists a regular version of the conditional 
probability given $\F_0$ on $\A$. We then use the notations of the introduction. 

\medskip

Recall that $UP_i=P_{i+1}U$ where $U$ is defined by $Uf=f\circ \theta$. For an adapted function $f\in L^2$ we thus have
$$
  f = \sum_{i=0}^\infty P_{-i}f +\E_{-\infty}(f)= \sum_{i=0}^\infty U^{-i}P_0U^if 
 + \E_{-\infty}(f)= \sum_{i=0}^\infty U^{-i}f_i +\E_{-\infty}(f)
$$
where $f_i = P_0U^if$, $i=0,1,\dots$. Therefore, since for 
every $n\ge 0$, $\E_0(\E_{-\infty}(f)) =\E_{-\infty}(f)$, we have
\begin{equation}\label{sum}
  \bar S_n(f) = S_n(f) - \E_0(S_n(f)) = \sum_{i=0}^{n-1}\sum_{j=1}^{n-i} U^jf_i.
\end{equation}

Denote, for $h\in L^1$, $Qh = \E_0(Uh)$. Then $Q$ is a Dunford-Schwartz operator (it is a contraction in all 
$L^p$,
$1\leq p\leq \infty$). Notice that $Q^nh=\E_0(U^nf)$. The use of the operator $Q$ is crucial in our proof. Its relevance 
to the problem is made more clear in the next section. 

\medskip

Let us recall several facts from ergodic theory that will be needed in the sequel.

\smallskip

By the Dunford-Schwartz (or Hopf) ergodic theorem (cf. \cite[Lemma 6.1]{Krengel}),  for
every $h\in L^1$, denoting $h^* = \sup_{n\geq 1}(1/n)\sum_{i=0}^{n-1} Q^i(|h|)$, we have
\begin{equation}\label{Hopf}
  \sup_{x>0} x\mu(h^*>x) \leq \|h\|_1.
\end{equation}

We will make use of the weak $L^2$-space $$L^{2,w}:=\{f\in L^1 ~:~\sup_{\lambda >0}\lambda^2 
\mu\{|f|\ge \lambda \}<\infty \}.$$
Recall that there exists a norm $\| \cdots \|_{2,w}$ on $L^{2,w}$ that makes it a Banach space and 
which is equivalent to the pseudo-norm $(\sup_{\lambda >0}\lambda^2
\mu\{|f|\ge \lambda \})^{1/2}$.

\smallskip

Then it follows from \eqref{Hopf}, that for every $h\in L^2$,
\begin{equation}\label{weakL2}
((h^2)^*)^{1/2}\in L^{2,w}.
\end{equation}

We obtain
\begin{lem}\label{lem}
Let $f$ be as above. We have
\begin{equation}\label{tightness}
(\E_0(\max_{1\le n\le N}\bar S_n^2(f)))^{1/2}\le \sqrt N \sum_{i=0}^\infty ((f_i^2)^*)^{1/2} <\infty\qquad \mbox{$\mu$-a.s.}
\end{equation}
In particular, if
$f$ satisfies the Hannan condition, then, by \eqref{weakL2} 
$$
\sup_{N\ge 1} \frac{\E_0(\max_{1\le n\le N}\bar S_n^2(f))}{ N} <\infty \qquad \mbox{$\mu$-a.s.}
$$
\end{lem}
{\bf Proof.} Let $N\ge n\ge 1$. From \eqref{sum} it follows that
$$
  |\bar S_n(f)|\le  \sum_{i=0}^{N-1} \max_{1\le k\le N} |\sum_{j=1}^{k} U^jf_i|.
$$
Notice that for every $i\geq 0$, the process $(U^jf_i)_j$ is a sequence 
of martingale increments.
We will use the Doob maximal inequality conditionally, in particular we will use
$$
  \big( E_0(\max_{n\leq N} |\bar S_n(f_i)|^2 ) \big)^{1/2} \leq
  2 \big[E_0(\bar S_N^2(f_i) )\big]^{1/2}.
$$
For $\mu$-a.e. $x\in X$ and every $i\ge 0$,
$(U^jf_i)_j$ remains a sequence of martingale increments under $\mu_x$. 
Denoting by $\|\cdot \|_{1,\mu_x}$ the norm in $L^2(\mu_x)$, it follows from the
Doob maximal inequality that 
\begin{gather*}
   \| \max_{n\leq N} |\bar S_n(f)|\|_{2,\mu_x}
  \leq \sum_{i=0}^{N-1} \| \max_{1\le n\le N} |\bar S_n(f_i)| \|_{2,\mu_x} \leq 2 \sum_{i=0}^{N-1} \|\bar S_N(f_i) \|_{2,\mu_x}
\end{gather*}
hence
\begin{gather*}
  \big( \E_0(\max_{n\leq N} |\bar S_n(f)|^2 ) \big)^{1/2} \leq
  2 \sum_{i=0}^{N-1} \big[\E(\bar S_N^2(f_i)) \big]^{1/2} =
  2 \sum_{i=0}^{N-1} \big[\E(\sum_{j=1}^N U^jf_i^2 ) \big]^{1/2} = \\
  = 2 \sum_{i=0}^{N-1} \big( \sum_{j=1}^N Q^jf_i^2 \big)^{1/2}
  \leq 2 \sqrt N \sum_{i=0}^\infty ((f_i^2)^*)^{1/2}.
\end{gather*}


Now, using \eqref{Hopf}, \eqref{weakL2}, we see that $\sum_{i=0}^\infty ((f_i^2)^*)^{1/2}$ is   in $L^{2,w}$, which finishes 
the proof.
\hfill $\square$

\bigskip




\noindent {\bf Proof of Theorem \ref{theo}.} By Hannan's condition $m=\sum_{k\ge 0}{P}_0(U^kf)$ is well
defined and $M_n=\sum_{k=1}^{n}U^km$ is a martingale with
stationary and ergodic increments.

Let $r\ge 1$. We have
$$
  f=\sum_{k=0}^r {P}_0(U^kf) -
\sum_{k=1}^r\big({\mathbb E}_0(U^kf)-{\mathbb E}_{-1}(U^{k-1}f)\big)+
{\mathbb E_{-1}(U^rf)}.
$$
Hence, denoting $m^{(r)}=\sum_{k=0}^r
{P}_0(U^kf)$ and $M_n^{(r)}=\sum_{l=1}^{n}U^lm^{(r)}$,
we obtain
\begin{equation*}
S_n-M_n= M^{(r)}_n-M_n -U^n(\sum_{k=1}^r{\mathbb
 E}_0(U^kf))+\sum_{k=1}^r{\mathbb E}_{0}(U^{k}f)
+\sum_{l=1}^{n} U^l({\mathbb E}_{-1}(U^rf))
\end{equation*}
and
\begin{gather}
\label{dec}S_n-M_n-{\mathbb E}_0(S_n) =\\
\nonumber = M^{(r)}_n-M_n -[U^n(\sum_{k=1}^r{\mathbb
 E}_0(U^kf))-{\mathbb E}_0(U^n(\sum_{k=1}^r{\mathbb
 E}_0(U^kf))] +
\\
\nonumber +\sum_{l=1}^{n} U^l({\mathbb E}_{-1}(U^rf))-{\mathbb E}_0(
\sum_{l=1}^{n} U^l({\mathbb E}_{-1}(U^rf)))
\end{gather}

By Doob maximal inequality, denoting $h^{(r)}:=(m-m^{(r)})^2$, we have
\begin{equation}\label{ineq0}
\E_0(\max_{1\le n \le N} (M^{(r)}_n-M_n)^2)\le 4 \sum_{1\le k \le N}Q^kh^{(r)}\le CN (h^{(r)})^*
\end{equation}
(recall that $h^* = \sup_{n\geq 1}(1/n)\sum_{i=0}^{n-1} Q^i(|h|)$).

Let $K>0$. Denote $Z^{(r)}=\sum_{k=1}^r{\mathbb E}_0(U^kf)$ and $Z^{(r)}_K=Z^{(r)}{\bf 1}_{|Z^{(r)}|>K}$
\begin{gather}
\nonumber \mathbb{E}_0(\max_{1\le n \le N}|U^n(\sum_{k=1}^r{\mathbb
 E}_0(U^kf)))-{\mathbb E}_0(U^n(\sum_{k=1}^r{\mathbb
 E}_0(U^kf))|^2)\\
 \nonumber \le 4 \E_0(\max_{1\le n \le N}|U^nZ^{(r)}|^2)
 \le 4K^2 + 4\E_0(\sum_{ n =1}^N|U^nZ^{(r)}_K|^2)\\
\label{ineq1} \le 4K^2 +4 \sum_{n=1}^N Q^n ((Z^{(r)}_K)^2) \le 4(K^2 + N ((Z^{(r)}_K)^2) ^*)
\end{gather}
To deal with the last term in \eqref{dec}, we apply  Lemma \ref{lem} to ${\mathbb E}_{-1}(U^rf)$, noticing that
in this case $f_i$ is replaced with ${P}_0 (U^i {\mathbb E}_{-1}(U^rf))={P}_0(U^{i+r}f)=f_{i+r}$
when $i\ge 1$ and for $i=0$, ${P}_0(\E_{-1}(U^rf))=0$). Hence
\begin{gather}\label{ineq2}
\E_0(\max_{1\le n\le N} |\sum_{l=0}^{n-1} U^l({\mathbb E}_{-1}(U^rf))-{\mathbb E}_0(
\sum_{l=0}^{n-1} U^l({\mathbb E}_{-1}(U^rf)))|^2)\le
N \sum_{i\ge r}((f_i^2)^*)^{1/2}.
\end{gather}
Combining \eqref{ineq0}, \eqref{ineq1} and \eqref{ineq2}, we obtain that for every $K>0$ and every $r\in \mathbb{N}$,
\begin{gather*}
\limsup_{N\to \infty}\frac{\E_0(\max_{1\le n\le N}(\bar S_{n}-M_n)^2)}{N}\le C(h^{(r)})^*+((Z^{(r)}_K)^2)^*
+ \sum_{i\ge r}((f_i^2)^*)^{1/2} \qquad \mbox{$\mu$-a.s.}
\end{gather*}

Now, $\|(((Z^{(r)}_K)^2)^*)^{1/2}\|_{2,w}\le C\|Z^{(r)}_K\|_2\underset{K\to \infty}
\longrightarrow 0$. Hence there exists a sequence $(K_l)$ going to infinity such that
$$
((Z^{(r)}_{K_l})^2)^*\underset{l\to \infty}
\longrightarrow 0 \qquad \mbox{$\mu$-a.s.}
$$
Hence
\begin{gather*}
\limsup_{N\to \infty}\frac{\E_0(\max_{1\le n\le N} (\bar S_{n}-M_n)^2)}{N}\le C(h^{(r)})^*+\sum_{i\ge r}((f_i^2)^*)^{1/2} <\infty \qquad \mbox{$\mu$-a.s.}
\end{gather*}
The second term clearly goes to 0 $\mu$-a.s., when $r\to \infty$ (by Lemma \ref{lem}), and the first one goes
to 0 $\mu$-a.s. (along a subsequence) by \eqref{weakL2}. \hfill $\square$

\medskip

\noindent {\bf Proof of Corollary \ref{cor}.} As noticed by Cuny-Peligrad \cite{CP}, the condition \eqref{SMW} implies
the Hannan condition and the fact that $\E_0(S_n)=o(\sqrt n)$ $\mu$-a.s., hence the result. \hfill $\square$

\section{Markov Chains}

In most of the literature, quenched limit theorems for stationary sequences use a Markov Chain setting:
the process is represented as a functional $(f(W_n))_n$ of a stationary and homogeneous Markov Chain $(W_n)$; the
limit theorem is said ``quenched'' it it remains true for almost every starting point.

Every (strictly) stationary sequence of random variables admits a Markov Chain representation. This has been observed
by Wu and Woodroofe in \cite{WW}, using an idea from \cite{Rosenblatt}. A remark-survey on equivalent representations of stationary
processes can be found in \cite{V2}. Here we show that the operator $Q$ introduced above leads to another Markov Chain
representation of stationary processes.


Let $(X,\mathcal{A},\mu)$ be a probability space and $\theta$ be an invertible bi-measurable
transformation of $X$ preserving the measure $\mu$.

\medskip


Let $\F\subset \A$ be a $\sigma$-algebra such that $\F\subset \theta^{-1}\F$. Denote $\E(\cdot |\F)$ the conditional expectation
with respect to $\F$ and define an operator $Q$ on
$L^\infty (X,\F,\mu)$ by \begin{equation}\label{operator}
Qh=\E(h\circ \theta |\F).
\end{equation}
 Then $Q$ is a positive contraction satisfying
$Q1=1$ and it is the dual of a positive contraction  $T$ of $L^1(X,\F,\mu)$, namely
$Tg= (\E(g|\F))\circ \theta^{-1}$.  By Neveu \cite[Proposition V.4.3]{Neveu}, if $X$ is a Polish space and
$\A$ the $\sigma$-algebra of its Borel sets, there exists a
transition probability $Q(x,dy)$ on $X\times \F$ such that for every $h\in L^{\infty}(X,\F,\mu)$,
$$
Qh=\int_Xh(y)Q(\cdot ,dy).
$$
Clearly, the transition probability $Q$ preserves the measure $\mu$, hence the canonical
Markov chain induced by $Q$, with initial distribution $\mu$, may be extended to $\mathbb{Z}$.

\medskip

Now define a sequence of random variables $(W_n)_{n\in \mathbb{Z}}$ defined from $(X,\A)$ to
$(X,\F)$ by  $W_n(x)=\theta^n(x)$. Then we have

\begin{prop}
Let $(X,\mathcal{A})$ be a Polish space with its Borel $\sigma$-algebra. Let $\mu$ be a probability on $\A$ and
$\theta$ be an invertible bi-measurable transformation of $X$ preserving the measure $\mu$. Let $\F\subset \A$ be a $\sigma$-algebra such that $\F\subset \theta^{-1}\F$. Then $(W_n)_{n\in \mathbb{Z}}$ is a Markov chain with state space $(X,\F)$, transition probability $Q$ (given by \eqref{operator}) and stationary distribution $\mu$. In particular, for every $f\in L^2(X,\F,\mu)$, the process $(f\circ \theta^n)$ is a functional of a stationary Markov chain.
\end{prop}
\noindent {\bf Proof.}\\
It suffices to show that for every $n\ge 1$, and any $\varphi_0,\ldots ,\varphi_n$ bounded measurable
functions from $\mathbb{R}$ to $\mathbb{R}$, we have
\begin{equation}\label{MarkProp}
\int_X \varphi_0 (W_0) \ldots \varphi_n(W_n)d\mu =\int_X\varphi_0(W_0)\ldots \varphi_{n-1}(W_{n-1})Q\varphi_n(W_{n-1})d\mu,
\end{equation}
the result for general blocks with possibly negative indices  follows by stationarity.\\
By definition of $(W_n)$, \eqref{MarkProp} holds for $n=1$.

\medskip

Assume that   \eqref{MarkProp} holds for a given $n\ge 1$. Let $\varphi_0,\ldots ,\varphi_{n+1}$ be bounded
$\F$-measurable
functions from $X$ to $\mathbb{R}$. Using the definition of $Q$, \eqref{MarkProp} for our given $n$, and stationarity,  we obtain
\begin{gather*}
\int_X \varphi_0 (W_0) \ldots \varphi_n(W_{n+1})d\mu=
\int_X\varphi_0 \varphi_1\circ \theta \ldots \varphi_{n+1}\circ
\theta^{n+1} d\mu\\
=\int_X\varphi_0(\theta^{-n})\ldots \varphi_{n-1}\circ \theta^{-1}\varphi_{n}\varphi_{n+1}\circ
\theta d\mu
=\int_X\varphi_0(\theta^{-n})\ldots \varphi_{n} Q\varphi_{n+1} d\mu \\=
\int_X\varphi_0\ldots \varphi_{n} \circ \theta^{n} Q\varphi_n\circ \theta^{n} d\mu
=\int_X\varphi_0(W_0)\ldots \varphi_{n} (W_n) Q\varphi_n(W_n) d\mu
\end{gather*}
which proves our result by induction. \hfill $\square$


\begin{thebibliography}{99}

\bibitem{B} P. Billingsley, Probability and Measure, (1995) Wiley.

\bibitem{CP} C. Cuny and M. Peligrad, \emph{Central limit theorem started
at a point for additive functionals of reversible Markov chains}, J. Theo. Prob.

\bibitem{DMV} J. Dedecker, F. Merlev\`ede and D. Voln\'y, \emph{On the weak invariance principle for non adapted sequences
under projective criteria},  J. Theo. Prob. 20, 2007, 971-1004.

\bibitem{D} O. Durieu, \emph{Independence of four projective criteria
for the weak invariance principle}, Al\'ea, Latin American Journal of
Probability and Mathematical Statistics 5, 2009, 21-27.

\bibitem{GP} M. Gordin and M. Peligrad, \emph{On the functional central limit theorem via
martingale approximation}, Bernoulli, 17, (1) 2011, 424-440.

\bibitem{H} E. J. Hannan, \emph{Central limit theorems for time series regression}, Z. Wahrscheinlichkeitstheorie
und Verw. Gebiete, 26, 1973, 157-170.

\bibitem{Krengel} U. Krengel, Ergodic theorems, de Gruyter Studies in Mathematics,
6. Walter de Gruyter \& Co., Berlin, 1985.

\bibitem{Neveu} J. Neveu, Bases math\'ematiques du calcul des probabilit\'es (French), 2nd ed.,  Masson et Cie,
Editeurs, Paris 1970 xii.

\bibitem{Rosenblatt} M. Rosenblatt, Markov Processes. Structure and asymptotic behaviour.
\emph{Springer-Verlag, New York-Heidelberg, 1971}.

\bibitem{VW1} D. Voln\'y and M. Woodroofe,  \emph{An example of non-quenched convergence in the conditional central limit theorem for
partial sums of a linear process}, Dependence in analysis, probability and number theory (The W. Philipp memorial volume),
Kendrick Press, 2010, 317-323

\bibitem{VW2} D. Voln\'y and M. Woodroofe,  \emph{Quenched central limit theorems for sums of stationary processes},
arXiv:1006.1795v1 [math.PR] 9 Jun 2010

\bibitem{V} D. Voln\'y,  \emph{On non ergodic versions of limit theorems}, Aplikace matematiky, 34, 1989, 351-363.

\bibitem{V2} D. Voln\'y,  \emph{Martingale approximation and optimality of some conditions for the central limit theorem},
J. Theo. Prob. 23, 2010, 888-903.

\bibitem{WW} W. Wu and M. Woodroofe, \emph{Martingale approximations for sums of stationary processes}, Ann.
Prob. 32 (2004), no. 2, 1674-1690.


\end{thebibliography}
\end{document}